\documentclass[reqno]{amsart}
\usepackage{lmodern,mathabx,mathtools}
\usepackage{tikz,shuffle}
\usepackage[inline]{enumitem}
\usepackage{tikz-cd}
\usepackage{cite,bookmark,cleveref}

\newtheorem{theorem}{Theorem}[section]
\newtheorem{proposition}[theorem]{Proposition}
\newtheorem{definition}[theorem]{Definition}
\newtheorem{lemma}[theorem]{Lemma}

\theoremstyle{remark}
\newtheorem{remark}[theorem]{Remark}

\DeclareMathOperator{\id}{id}
\def\WickOp#1{:\mkern-4mu#1\mkern-4mu:\mskip-4mu plus 2mu}
\def\Wick#1{:\mkern-4mu#1\mkern-4mu:\mkern6mu}
\def\1{\mathbf 1}

\begin{document}

%%%%%%%%%%%%%%%%%%%%%%%%%%%%%%%%%
%%%%%%%%%%%%%%%%%%%%%%%%%%%%%%%%%

\title[Wick polynomials in non-commutative probability]{Wick polynomials in non-commutative probability.\\
A group-theoretical approach}

%%%%%%%%%%%%%%%%%%%%%%%%%%%%%%%%%
%%%%%%%%%%%%%%%%%%%%%%%%%%%%%%%%%

\author[K.~Ebrahimi-Fard]{K.~Ebrahimi-Fard}
\address{Department of Mathematical Sciences, NTNU, Trondheim, Norway.}
\email{kurusch.ebrahimi-fard@ntnu.no}
\urladdr{https://folk.ntnu.no/kurusche/}

\author[F.~Patras]{F.~Patras}
\address{Univ.~C\^ote d'Azur, CNRS, UMR 7351, Parc Valrose, Nice, France.}
\email{patras@unice.fr}
\urladdr{www-math.unice.fr/$\sim$patras}

\author[N.~Tapia]{N.~Tapia}
\address{Weierstrass Institute and Technische Universit\"at Berlin, Berlin, Germany.}
\email{tapia@wias-berlin.de}
\urladdr{http://wias-berlin.de/people/tapia/}

\author[L.~Zambotti]{L.~Zambotti}
\address{LPSM, Sorbonne Universit\'e, Sorbonne Université, Université de Paris, CNRS, Paris, France}
\email{lorenzo.zambotti@upmc.fr}
\urladdr{http://www.lpma-paris.fr/pageperso/zambotti}

\thanks{N.T.~kindly acknowledges support from the European Research Council for Informatics and Mathematics through contract ERCIM 2018-10.}
\keywords{Wick polynomials, monotone cumulants, free cumulants, boolean cumulants,
formal power series, combinatorial Hopf algebra, shuffle algebra, group actions}
\subjclass[2020]{16T05, 16T10, 16T30, 17A30, 46L53, 46L54}

\maketitle

\section{Introduction}
\label{sect:intro}

Moment-cumulant relations and Wick products play a central role in probability theory and related fields \cite{Akhiezer65,peccati2011wiener}. In classical probability, cumulant sequences $(c_n)_{n\in{\mathbf N}^\ast}$ linearize the notion of independence of random variables: if two random variables, $X,Y$, with moments of all orders are independent then for $n\geq 1$, $c_n(X+Y)=c_n(X)+c_n(Y)$. Wick polynomials, Wick products and chaos expansions are related to cumulants. Indeed, recall for example that given a random variable $X$ with moments of all orders, the Wick polynomial $W(X^n)$ is the coefficient of $\frac{t^n}{n!}$ in the expansion of $\exp(tX-K(t))$, where $K(t)$ is the exponential generating series of cumulants.

Voiculescu's theory of free probability \cite{Voi97,VDN92} provides the paradigm of a non-commutative probability theory, where the notion of freeness replaces the classical concept of probabilistic independence. Speicher showed that free cumulants linearise Voiculescu's notion of freeness. See \cite{MS17,NS06} for detailed introductions. Following Voiculescu's ideas, various authors \cite{HS11,M2001,Spe97,B1986,SW97} considered different types of independences (Boolean, monotone, and others), each characterized by particular moment-cumulant relations with explicit combinatorial descriptions given in terms of different types of set partitions. Relations between the different brands of cumulants were thoroughly explored by Arizmendi et al.~in \cite{AHL+15}. Free and Boolean Wick polynomials have been introduced in this setting by Anshelevich \cite{Ans04,Ans09a,Ans09}.

In a previous paper \cite{EPTZ18}, the authors presented a Hopf-algebraic framework describing both the combinatorial structure of the classical moment-cumulant relations as well as the related notions of Wick polynomials and Wick products. The approach is based on convolution products of linear functionals defined on a coalgebra and encompasses the multi-dimensional extension of the moment-cumulant relations. In this framework, classical Wick polynomials result from a Hopf algebra deformation under the action of linear automorphisms induced by multivariate moments associated to an arbitrary family of random variables with moments of all orders.

In a series of recent papers \cite{EP15,EP17,EP18}, two of us explored relations between multivariate moments and free, boolean and monotone cumulants as well as relations among latter in non-commutative probability theory by studying a particular graded connected Hopf algebra $H$ defined on the double tensor algebra over a non-commutative probability space $(A,\varphi)$. In this approach, the associated set partitions (non-crossing, interval and monotone, respectively) appear through the evaluation of elements of the  group $G$ (Lie algebra $\mathfrak g$) of (infinitesimal) Hopf algebra characters on words. %It turns out that the Hopf algebra $H$ describes not only the moment-cumulant relations in non-commutative probability, but more generally, it encodes an action of the group of formal diffeomorphisms on the group of invertible power series in non-commuting variables. This relation will be further explored in a forthcoming paper.

In the paper at hand we revisit from a Hopf theoretic point of view the theory of free, boolean and conditionally free Wick polynomials. The relevance of shuffle group actions and structures in the sense of \cite{EP17} is also emphasized. 

The article is organized as follows. In \Cref{sect:classical}, we recall the definitions of
classical cumulants and Wick polynomials. In \Cref{sect:cumulants} we do the same for free and
boolean cumulants.
\Cref{sect:nc-Wick} defines free Wick polynomials using the Hopf algebraic approach. The new definition is shown to extend Anshelevich's definition of multivariate free Appell polynomials.
At the beginning of \Cref{sect:nsShuffle} we introduce the shuffle-theoretic framework allowing to
deal with non-commutative moment-cumulant relations and the corresponding non-commutative Wick polynomials. \Cref{sect:ncrev} revisits accordingly moment-cumulant relations in non-commutative probability theory following mainly the references \cite{EP15,EP18a,EP18}. \Cref{sect:shW} develops shuffle calculus for free Wick polynomials.
In \Cref{sec:bwickpoly} boolean Wick polynomials are also introduced and analysed from this point of
view. \Cref{sect:cfreeW} uses the same approach to define conditionally free Wick polynomials. In
\Cref{sect:group} we show how the three notions of non-commutative Wick polynomials can be related
through comodule structures and the induced group actions. \Cref{sect:products} shows how the classical notion of Wick products generalises naturally to the non-commutative setting, inducing three new associative algebra structures on the tensor algebra over a non-commutative probability space.
Finally, in \Cref{ssect:tensorc}, we show using a Hopf algebraic approach how the definition of classical cumulants lifts to the notion of tensor cumulants for random variables in a non-commutative probability space. In the following \Cref{sect:Wick}, we explain how this leads to the definition of tensor Wick polynomials. These two sections extend the results of \cite{EPTZ18} from the classical to the tensor framework.

Below, $\mathbb{K}$ denotes the base field of characteristic zero over which all algebraic structures are defined. All (co-)algebras are (co-)associative and (co-)unital unless otherwise stated.

%%%%%%%%%%%%%%%%%%%%%%%%%%%%%%%%%
%%%%%%%%%%%%%%%%%%%%%%%%%%%%%%%%%

\section{Cumulants and Wick polynomials}
\label{sect:classical}

Let us first recall briefly the definition of classical cumulants and Wick polynomials. Let $X$ be a real-valued random variable, defined on a probability space $(\Omega,\mathcal F,\mathbb P)$, with finite moments of all orders, i.e., such that $m_n \coloneq \mathbb EX^n<\infty$ for all $n>0$. Its exponential moment-generating function is defined as a power series in $t$
\begin{equation}
  M(t)\coloneq\mathbb E\, \mathrm \exp({tX})
    = 1+ \sum_{n > 0}m_n\frac{t^n}{n!}.
  \label{eq:momgen}
\end{equation}
If we assume suitable growth conditions on the coefficients $m_n$ so that the above series has a positive radius of convergence, then this power series defines a function of class $C^\infty$ around the origin, and the moments $m_n$ can be recovered from it by differentiation.

The exponential cumulant-generating function
\begin{equation*}
  K(t)\coloneq \sum_{n>0} c_n \frac{t^n}{n!}
\end{equation*}
is a power series in $t$ defined through the classical exponential relation between moments and cumulants
\begin{equation}
  M(t)=\exp\left( K(t) \right).
\label{eq:ccdef}
\end{equation}
Using standard power series manipulations, this equation rewrites:
\begin{equation}
  m_n=\sum_{\pi\in P(n)}\prod_{B\in\pi}c_{|B|}.
\label{eq:cmcrel}
\end{equation}
Here, $P(n)$ denotes the collection of all set partitions, $\pi\coloneq\{B_1,\ldots,B_l\}$, of the set $[n]\coloneq\{1,\dotsc,n\}$, where the block $B_i \in \pi$ contains $|B_i|$ elements. In general, for a finite subset $U \subset \mathbb N$ we denote by $P(U)$ the collection of all set partitions of $U$.

Let $(X_1,\dotsc,X_p)$ be a finite collection of real-valued random variables defined on a common probability space, such that all the  moments $m_{\mathbf n}\coloneq\mathbb E[X_1^{n_1}\dotsm X_p^{n_p}]$, where $\mathbf n\coloneq(n_1,\dotsc,n_p)\in\mathbb N^p$ is a multi-index, exist.
We may consider a multivariate extension of \eqref{eq:momgen}, namely
\begin{align}
    M(t_1,\dotsc,t_p)
  &\coloneq\mathbb E\, \mathrm \exp({t_1X_1+\dotsb+t_pX_p}) \nonumber \\
        &\eqcolon\sum_{\mathbf n}m_{\mathbf n}\frac{t^{\mathbf n}}{\mathbf n!}, \label{eq:mvmgen}
\end{align}
where  $t^{\mathbf n}\coloneq t_1^{n_1}\dotsm t_p^{n_p}$ and $\mathbf n!\coloneq n_1!\dotsm n_p!$. As before, the cumulant-generating function is defined by a relation analogous to \eqref{eq:ccdef}, and its coefficients are related to the moments in a way analogous to \eqref{eq:cmcrel}. This relation will be revisited in the following sections.

There exists a particular family of polynomials associated to a random variable $X$ with finite moments of all orders, called \emph{Wick polynomials} and denoted here by $W_n(x)$, $n\ge 0$. It turns out to be the unique family of polynomials such that $W_0(x)=1$ and
\[
  \mathbb E\, W_n(X)=0, \qquad
  \frac{\mathrm d}{\mathrm dx}W_n(x)=nW_{n-1}(x),
\]
for all $n>0$.
The latter defining property means that $(W_n)_{n\ge0}$ qualifies as a sequence of Appell polynomials \cite{Appell1880}. For example, if $X$ is a standard Gaussian random variable, this family coincides with the Hermite polynomials. These polynomials are interesting for physics. In particular, the \emph{Wick exponential}
\[
  \WickOp{\exp}(tX)\coloneq\sum_{n\ge0}W_n(X)\frac{t^n}{n!}
  =\frac{\exp({tX})}{\mathbb E \exp({tX})}=\exp({tX-K(t)})
\]
is closely related to moment- and cumulant-generating functions. In fact, this relation can be used to define Wick polynomials since the exponential power series in $t$ serves as a generating function. %Thus, the above identity relates the generating function of Wick polynomials to the moment- and cumulant-generating functions. We also note that both \eqref{eq:cWickcenter,eq:cWickdiff} can be deduced if we use \eqref{eq:cWickgen} as a starting point instead.

The polynomial
$$
  \Wick{X^n}\coloneq W_n(X)
$$
is called the $n$-th \emph{Wick power} of $X$. For example,
\[
  \Wick{X}=X-\mathbb EX,
  \quad
  \Wick{X^2}=X^2-2X\,\mathbb EX+2(\mathbb EX)^2-\mathbb EX^2,
  \quad
  \dotsc
\]
In general, these explicit expansions can be recursively obtained from the change of basis relation
\begin{equation}
\label{eq:wickbasis}
  x^n=\sum_{j=0}^n\binom{n}{j}W_j(x)\,m_{n-j}.
\end{equation}
The latter can be generalized to finite collections $(X_1,\dotsc,X_p)$ of random variables in a way analogous to \eqref{eq:mvmgen}.

%%%%%%%%%%%%%%%%%%%%%%%%%%%%%%%%%
%%%%%%%%%%%%%%%%%%%%%%%%%%%%%%%%%

\section{Free and boolean cumulants}
\label{sect:cumulants}

Voiculescu introduced free probability theory in the 1980s \cite{Voi97,VDN92}\footnote{The referee pointed us to the early reference \cite{A1982} for construction of a free product state on the free product of a family of $C^\ast$-algebras.}. In this theory the classical notion of independence is replaced by the algebraic notion of freeness.  A family of unital subalgebras $(B_i:i\in I)$ of a non-commutative probability space \((A,\varphi)\) is called \emph{freely independent} (or free), if $\varphi(a_1 \cdot_{\!\scriptscriptstyle{A}} \dotsm \cdot_{\!\scriptscriptstyle{A}}  a_n)=0$ whenever $\varphi(a_j)=0$ for all $j=1,\dotsc,n$ and $a_j\in B_{i_j}$ for some indices $i_1\neq i_2 \neq\dotsb\neq i_n$.

Speicher introduced the notion of \emph{free cumulants} \cite{Spe97} as the right analogue of the classical cumulants in the theory of free probability, allowing for a more tractable characterisation of Voiculescu's notion of freeness. Free cumulants are defined by a formula analogous to \eqref{eq:cmcrel} where the lattice $P$ of set partitions is replaced by the lattice $\operatorname{NC}$ of non-crossing partitions:
\begin{equation}
  \varphi(a_1 \cdot_{\!\scriptscriptstyle{A}} \dotsm \cdot_{\!\scriptscriptstyle{A}} a_n)
  =\sum_{\pi\in \operatorname{NC}([n])}\prod_{B\in\pi}k(a_B).
  \label{eq:fmcrel}
\end{equation}
As above, we set $k (a_B)\coloneq k (a_{i_1}, \ldots, a_{i_{|B|}})$, for $B=\{i_1 < \cdots < i_{|B|}\}$, to be the multivariate free cumulant of order $|B|$. Free cumulants reflect freeness in the sense that they vanish whenever the involved random variables belong to different freely independent subalgebras.

Relation \eqref{eq:fmcrel} between  moments and free cumulants can be concisely expressed in terms of their ordinary generating functions. Indeed, given $a_1,\dots,a_n$ in $A$, introduce non-commuting variables $w_1,w_2,\dots,w_n$ and the generating functions
\[
  M(w)\coloneq1+\sum_{\mathbf n}\varphi(a_{\mathbf n})\,w_{\mathbf n},
  \quad
  R(w)\coloneq\sum_{\mathbf n} k(a_{\mathbf n})\,w_{\mathbf n}.
\]
Here we define $\varphi(a_{\mathbf n})\coloneq\varphi(a_{n_1} \cdot_{\!\scriptscriptstyle{A}} \dotsm \cdot_{\!\scriptscriptstyle{A}} a_{n_p})$ for the multi-index $\mathbf n \coloneq (n_1,\dotsc,n_p)\in [n]^p,\ p\in{\mathbb N}^\ast$, and similarly for $k(a_{\mathbf n})=k(a_{n_1}, \ldots, a_{n_p})$.
Then, \eqref{eq:fmcrel} is summarized by the intriguing identity \cite{Ans04,NS06}
$$
  M(w)=1+R(z),
$$
where the substitution
\begin{equation}
\label{substitut}
  z_i\coloneq w_iM(w)
\end{equation}
is in place on the righthand side.

The fact that the random variables under consideration do not commute entails that we are able to consider several other notions of independence in addition to Voiculescu's freeness. For example, the notion of \emph{boolean cumulants} appears naturally in the context of the study of stochastic  differential equations \cite{vW1973}. Speicher and Woroudi  \cite{SW97} defined the multivariate boolean cumulants, $ b(a_1, \ldots , a_n)$, and the corresponding relations with moments in the context of non-commutative probability theory in terms of the following recursion
\[
  \varphi(a_1 \cdot_{\!\scriptscriptstyle{A}} \dotsm \cdot_{\!\scriptscriptstyle{A}} a_n)
  =\sum_{j=1}^n b(a_1, \ldots , a_j) \, \varphi(a_{j+1} \cdot_{\!\scriptscriptstyle{A}} \dotsm\cdot_{\!\scriptscriptstyle{A}} a_n).
\]
While the combinatorics of free cumulants is described by the lattice of non-crossing partitions, the relation between moments and boolean cumulants can be expressed by using the lattice $\operatorname{Int}$ of interval partitions:
\[
  \varphi(a_1 \cdot_{\!\scriptscriptstyle{A}} \dotsm \cdot_{\!\scriptscriptstyle{A}} a_n)
  =\sum_{\pi\in \operatorname{Int}([n])}\prod_{B\in\pi} b(a_B),
  \qquad b(a_B)\coloneq b(a_{i_1}, \ldots, a_{i_{|B|}}).
\]
Using the multi-index notation from above, these relations can be encapsulated in a single identity by introducing the generating function
\[
  \eta(w)\coloneq\sum_{\mathbf n}b(a_{\mathbf n})w_{\mathbf n},
\]
yielding the simple expression \cite{Ans09a, NS06}
$$
  M(w)=1+\eta(w)M(w).
$$

Observe that in this case, as opposed to the functional equation describing the relation between moments and free cumulants, there is no substitution such as \eqref{substitut} to be made.

% Finally, we consider a third type of independence known as \emph{monotone independence} \cite{AHL+15,HS11}. It is also captured through a set of cumulants, $r(a_B)=r(a_{i_1}, \ldots , a_{i_{|B|}})$, for $B=\{i_1 < \cdots < i_{|B|}\}$, defined in terms of  moments by the relations
% \[
%   \varphi(a_1\cdot_{\!\scriptscriptstyle{A}} \dotsm\cdot_{\!\scriptscriptstyle{A}} a_n)
%   =\sum_{\pi\in \operatorname{NC}([n])}\frac{1}{\tau(\pi)!}\prod_{B\in\pi}r(a_B).
%   \]
% Here $\tau(\pi)$ is the forest of rooted trees encoding the nesting structure of the non-crossing partition $\pi$ and $\tau(\pi)!$ is the corresponding tree factorial. See \cite{AHL+15} for details. In this case, the description in terms of generating functions is more involved and requires solving a certain ODE \cite[Proposition 4.7]{HS11}. We mention that in the shuffle algebra context this can be simplified to a linear initial value problem (see Section \ref{sect:ncrev} further below).

Surprisingly, the relation between  moments and the different types of cumulants can be described concisely as the action of linear maps on the double tensor algebra. For this, two of us introduced in \cite{EP15}, a different coproduct which allows to express these relations in a way similar to the presentation of the preceding sections.

%%%%%%%%%%%%%%%%%%%%%%%%%%%%%%%%%
%%%%%%%%%%%%%%%%%%%%%%%%%%%%%%%%%

\section{Free Wick polynomials}
\label{sect:nc-Wick}
In \cite{EP15,EP18,EPTZ18} an approach in terms of Hopf algebras to the moment-cumulant relations in
both classical and non-commutative probability was introduced. It permits to describe
moment-cumulant relations in a rather different way, avoiding the use of generating functions.
\begin{definition}
  A non-commutative probability space $(A,\varphi)$ consists of a unital algebra $A$ together with a unital map $\varphi\colon A \to \mathbb{K}$, i.e., $\varphi(1_A)=1$.
  \label{dfn:ncproba}
\end{definition}

To avoid ambiguities we also denote the product of elements $a,b$ in the algebra $A$ by $m_{A}(a \otimes b)\eqcolon a \cdot_{\!\scriptscriptstyle{A}} b$.
We still write $m_A$ for the iterates
$$
  m_A\colon a_1\otimes\dots \otimes a_n\longmapsto a_1 \cdot_{\!\scriptscriptstyle{A}} \dots \cdot_{\!\scriptscriptstyle{A}} a_n.
$$

%In analogy with the previous paragraphs, it is useful to think of $A$ as the algebra of random variables with finite moments of all orders defined on a common (measure-theoretical) probability space, and of $\varphi$ as the associated expectation operator. An important remark is in order. Even though our aim is --at this point-- to capture the moment-cumulant relation in classical, that is, commutative probability,
Notice that we do not require the algebra $A$ to be commutative. The elements of $A$ should be thought of in general as non-commutative random variables and the map $\varphi$ plays then the role of the expectation map.
Elements in $A$ can represent, for example, operator-valued random variables such as those appearing in the Fock space approach to Quantum Field Theory \cite{EP03}.

We consider the non-unital tensor algebra over $A$
\[
  T(A)\coloneq\bigoplus_{n>0} A^{\otimes n}
\]
and we denote elements of $T(A)$ using word notation ($a_1\cdots a_n=a_1\otimes\dots \otimes a_n$). It is graded by the number of letters, i.e., the length of a word. The unitalization of $T(A)$ follows from adding the empty word $\1$ and is denoted by $\overline T(A)=T_0(A)\oplus T(A)\coloneq \mathbb K\1\oplus T(A)$.
The product on $T(A)$ (resp.~$\overline T(A)$) is given by concatenation of words, $\mathrm{conc}(w_1 \otimes w_2)\coloneq w_1w_2$, for $w_1,w_2 \in T(A)$ (with the empty word $\1$ being the unit).
Let $A$ be an algebra and consider the double tensor algebra $\overline T(T(A))$ over $A$.
On \(\overline T(T(A))\) we also consider the concatenation product, but we denote it with a vertical bar in order to distinguish it from concatenation in \(T(A)\), i.e., \(\mathrm{conc}(w_1\otimes w_2)=w_1|w_2\) for \(w_1,w_2\in\overline T(T(A))\).

Given a subset $U\subset\mathbb N$, \emph{an interval or connected component} of $U$ is a maximal sequence of successive elements in $U$. For a subset $S\subseteq[n]$ we denote by $J_1^S,\dotsc,J_{k(S)}^S$ the connected components of $[n]\setminus S$, ordered in increasing order of their minimal element. For notational convenience, we will often omit making explicit the dependency on $S$ of the number of these connected components and, when there is no risk of confusion, will write simply $J_1^S,\dotsc,J_k^S$ for $J_1^S,\dotsc,J_{k(S)}^S$.

\begin{definition}\label{def:coprod}
  The map $\Delta\colon T(A)\to \overline{T}(A)\otimes\overline T(T(A))$ is defined by
\begin{equation}
\label{eq:theCoprod}
  \Delta(a_1\cdots a_n)
  \coloneq a_1\cdots a_n \otimes \1 
  			+ \1 \otimes a_1\cdots a_n 
	+ \sum_{\substack{S\subsetneq[n]\\S \neq \emptyset}}
				a_S\otimes a_{J^S_1}\vert\dotsm\vert a_{J^S_k}.
\end{equation}
It has a unique multiplicative extension $\Delta\colon\overline T(T(A))\to\overline T(T(A))\otimes\overline T(T(A))$ such that $\Delta(\1)=\1\otimes\1$.
\end{definition}

Note that in the sum on the righthand side of \eqref{eq:theCoprod}, we have inserted the concatenation product in \(\overline T(T(A))\) between the words corresponding to the connected components ${J^S_1},\dotsc,{J^S_k}$ associated to the non-empty set $S\subsetneq[n]$, that is, whereas $a_S \in {T}(A)$ we have $a_{J^S_1}\vert\dotsm\vert a_{J^S_k} \in T(T(A))$.

\begin{theorem}[\cite{EP15}]
  The unital double tensor algebra $\overline T(T(A))$ equipped with $\Delta$ is a non-commutative non-cocommutative connected graded Hopf algebra.
\end{theorem}

Extending our approach to classical Wick polynomials into the non-commutative realm, we introduce an endomorphism of the double tensor algebra $\overline T(T(A))$. This provides, among others, a new way of introducing the non-commutative Wick (a.k.a.~free Appell) polynomials appearing in the work of Anshelevich \cite{Ans04}, as explained below. 

Suppose that $(A,\varphi)$ is a probability space. We define the map $\Phi \colon \overline T(T(A))\to \mathbb{K}$ as the unique unital multiplicative extension of the linear map $\phi$ defined on $T(A)$ by $\phi(a_1\dotsm a_n) \coloneq \varphi(a_1\cdot_{\!\scriptscriptstyle{A}} \dotsm\cdot_{\!\scriptscriptstyle{A}} a_n)$. Since $\Phi$ is --by definition-- a Hopf algebra character, it is an invertible element in the corresponding convolution algebra. Its convolution inverse, denoted $\Phi^{-1}$, is the unique character on the double tensor algebra such that $\Phi^{-1}*\Phi=\Phi*\Phi^{-1}=\varepsilon$. 
Here, \(\varepsilon\colon\overline T(T(A))\to\mathbb K\) denotes the \emph{counit}, defined as the unique multiplicative map such that \(\ker\varepsilon=T(T(A))\), and which acts as the neutral element for the convolution product.
In other words, the map \(\varepsilon\) is such that \(\varepsilon(\mathbf 1)=1\) and vanishes otherwise, and \(\varepsilon(w_1|w_2)=\varepsilon(w_1)\varepsilon(w_2)\) for all \(w_1,w_2\in\overline T(T(A))\).

\begin{definition}
  \label{dfn:fWick}
  The free Wick map $\mathrm{W} \colon\overline T(T(A))\to\overline T(T(A))$ is defined by
  \[
    \mathrm{W} \coloneq({\id}\otimes\Phi^{-1})\Delta,
  \]
  or, implicitly, by
  \[
    \id =(\mathrm{W}\otimes\Phi)\Delta.
  \]
\end{definition}
 We call \emph{free Wick polynomials} the family $\{\mathrm{W}(a_1\cdots a_n)$, $a_i\in A$, $i=1,\ldots, n\}$.

\begin{proposition}
The free Wick map is multiplicative, i.e., for words $w,w'\in T(A)$,
$$
  \mathrm{W}(w|w')=\mathrm{W}(w)|\mathrm{W}(w').
$$
We recall that \(a|b\) denotes the concatenation of \(a\) and \(b\) in \(\overline T(T(A))\).
\end{proposition}

\begin{proof}
As the identity map $\id$ and $\Phi^{-1}$ are both multiplicative, using Sweedler's notation, $\Delta(w)= \sum w^{(1)}\otimes w^{(2)}$, for the coproduct defined in \eqref{eq:theCoprod}:
\begin{align*}
  \mathrm{W}(w|w')
    &=(\id\otimes \Phi^{-1})\Delta(w|w')\\
    &=(\id\otimes \Phi^{-1})(\Delta(w)\Delta(w'))\\
    &=\sum \sum (w^{(1)}|{w'}^{(1)})\,\Phi^{-1}(w^{(2)})\,\Phi^{-1}({w'}^{(2)})
    =\mathrm{W}(w)|\mathrm{W}(w').
\end{align*}
\end{proof}

The compositional inverse of $\mathrm{W}$, denoted $\mathrm{W}^{\circ -1}$, is given by
\[
  \mathrm{W}^{\circ -1}=({\id}\otimes\Phi)\Delta.
\]

From Definition \ref{dfn:fWick}, we also obtain that the usual monomials in $\overline T(A)$ can be expressed in terms of free Wick polynomials:
\begin{equation}\label{eq:ans04}
  a_1\dotsm a_n=\sum_{S\subseteq[n]}  \mathrm{W}(a_S)\Phi(a_{J_1^S})\dotsm\Phi(a_{J_k^S}).
\end{equation}
Note that $\mathrm{W}$ restricts to an automorphism of $\overline T(A)$. By \cite[Prop. 3.12]{Ans04}, our free Wick polynomials agree with Anshelevich's free Appell polynomials since our formula \eqref{eq:ans04} coincides with formula \cite[Formula (3.42)]{Ans04}.

Here are some low-degree computations:
\begin{align}
  \mathrm{W}(a_1) &= a_1-\varphi(a_1)\1, \nonumber\\
  \mathrm{W}(a_1a_2) &= a_1a_2-\varphi(a_2)a_1-\varphi(a_1)a_2
        - \big(\varphi(a_1\cdot_{\!\scriptscriptstyle{A}} a_2)-2\varphi(a_1)\varphi(a_2)\big)\1,  \nonumber\\
  \begin{split}
  \mathrm{W}(a_1a_2a_3) &= a_1a_2a_3-\varphi(a_3)a_1a_2-\varphi(a_2)a_1a_3-\varphi(a_1)a_2a_3\\
  &\quad- \big(\varphi(a_2\cdot_{\!\scriptscriptstyle{A}} a_3)
    - 2\varphi(a_2)\varphi(a_3)\big)a_1
      + \varphi(a_1)\varphi(a_3)a_2
    - \big(\varphi(a_1\cdot_{\!\scriptscriptstyle{A}} a_2)\\
  &\quad -2\varphi(a_1)\varphi(a_2)\big)a_3
    - \big(\varphi(a_1\cdot_{\!\scriptscriptstyle{A}} a_2\cdot_{\!\scriptscriptstyle{A}} a_3)
    - 2\varphi(a_1)\varphi(a_2\cdot_{\!\scriptscriptstyle{A}} a_3)\\
  &\quad -2\varphi(a_3)\varphi(a_1\cdot_{\!\scriptscriptstyle{A}} a_2)
    - \varphi(a_2)\varphi(a_1\cdot_{\!\scriptscriptstyle{A}} a_3)
  + 5\varphi(a_1)\varphi(a_2)\varphi(a_3)\big)\1.
  \end{split} \label{freeWick3}
\end{align}
The computation of the third order polynomial \eqref{freeWick3} is somewhat subtle and should be compared with the expression \eqref{Wick3} below.

The free Wick polynomials inherit immediately from their Hopf algebraic definition a key property of classical Wick polynomials.

\begin{lemma}
The Wick polynomials $\mathrm{W}$ in Definition \ref{dfn:fWick} are centred. That is,
\[
  \Phi \circ \mathrm{W}
  =(\Phi\otimes\Phi^{-1})\Delta
  =\Phi*\Phi^{-1}=\varepsilon.
\]
\end{lemma}

\begin{definition}
Let us call universal polynomial $P=P(x_1,\ldots,x_n;\gamma)$ for non-commutative probability spaces any linear combination of symbols
$$
  \gamma(X^{\bullet}_{J_1}) \cdots \gamma(X^\bullet_{J_p})X_I,
$$
where $I\coprod J_1\coprod \cdots \coprod J_p$ is a partition of $[n]$ and $\gamma$ takes values in $\mathbb{K}$.
\end{definition}

To a universal polynomial $P$ together with a non-commutative probability space $(A,\varphi)$ and elements $a_1,\ldots,a_n \in A$, we associate the element $P(a_1,\ldots,a_n;\varphi)\in\overline{T}(A)$  obtained from $P$ by replacing $X_I$ with the tensor monomial $a_{i_1}\cdots a_{i_k}$, where $I=\{i_1,\ldots,i_k\}$, and $X_J^\bullet$ with $a_{j_1}\cdot_{\!\scriptscriptstyle{A}} \cdots \cdot_{\!\scriptscriptstyle{A}} a_{j_l}$, where $J=\{j_1,\ldots,j_l\}$.

A family $(f_{(A,\varphi)})$ of linear endomorphisms of $\overline{T}(A)$, where $(A,\varphi)$ runs over non-commutative probability spaces, is called \emph{universal} if its action on words $a_1\cdots a_n$ is given by universal polynomials. The Wick map, $\mathrm{W}$, the inverse Wick map, $\mathrm{W}^{\circ -1}$, the moment map, and the cumulant maps are examples of universal families.

Now, given $(A,\varphi)$, we define a formal derivation with respect to an element $a \in A$ as follows. Fix a decomposition
$A=\mathbb{K}a \oplus A^\prime$. Denote by $\zeta_a \colon T(A)\to \mathbb{K}$ the linear map defined by $\zeta_a(a)\coloneq1$, $\zeta_a(b)\coloneq0$ for $b\in A^\prime$ and $\zeta_a(w)\coloneq0$ for every word $w=a_1\cdots a_n$, $a_i \in A$, $n\geq 2$. This map (which depends on the chosen direct sum decomposition of $A$) is then extended as an infinitesimal character to the double tensor algebra. We set
$$
  \partial_a\colon\overline T(T(A))\to\overline T(T(A)), 
  \qquad \partial_a\coloneq(\zeta_a\otimes\id)\Delta.
$$
Observe that for any word $w=a_1\cdots a_n \in T(A)$ where $a_j=a$ or $a_j \in A^\prime$, we then get
\[
  \partial_a(w)=\sum_{j:a_j=a} a_1\dotsm a_{j-1}\vert a_{j+1}\dotsm a_n.
\]
For example, if $w,w_1,w_2\in T(A^\prime)$, then
\[
  \partial_a(aw)=w=\partial_a(wa),\quad \partial_a(w_1aw_2)
  =w_1\vert w_2,\quad\partial_a(awa)=aw+wa.
\]
Since $\zeta_a$ is infinitesimal, $\partial_a$ turns out to be a derivation on $\overline T(T(A))$.

\begin{theorem}
  The Wick map $\mathrm{W}$ is the unique family of algebra automorphisms of $\overline T(T(A))$, where $(A,\varphi)$ runs over non-commutative probability spaces, such that
  \begin{itemize}
  \item The restrictions of $\mathrm{W}$ to $\overline{T}(A)$ form a universal family.
  \item The map $\mathrm{W}$ is centred, $\Phi \circ \mathrm{W}=\varepsilon$, with $W(\1)=1$ in particular.
  \item For any $a \in A$ and any direct sum decomposition $A={\mathbb K} a \oplus A^\prime$
$$
    \partial_a\circ \mathrm{W}=\mathrm{W}\circ\partial_a.
$$
  \end{itemize}
\end{theorem}

\begin{proof}
  The first two statements were already shown.
  The third one follows from the coassociativity of the coproduct:
  \begin{align*}
    \partial_a\circ \mathrm{W}&=(\partial_a\otimes\Phi^{-1})\Delta\\
    &=(\zeta_a\otimes{\id}\otimes\Phi^{-1})(\Delta\otimes{\id})\Delta\\
    &=(\zeta_a\otimes{\id}\otimes\Phi^{-1})({\id}\otimes\Delta)\Delta\\
    &=\mathrm{W}\circ\partial_a.
  \end{align*}
Uniqueness follows from the fact that these three properties define the universal family $\mathrm{W}$ by induction. Given an integer $n$, choose for example a family $a_1,\ldots,a_n$ of linearly independent free random variables in a non-commutative probability space $(A,\varphi)$. Use then an adapted direct sum decomposition $A={\mathbb K}a_1\oplus\dots\oplus {\mathbb K}a_n\oplus A^{\prime\prime}$ to define the derivations. The knowledge of the
$$
  \partial_{a_i}\mathrm{W}(a_1\cdots a_n)
  =\mathrm{W}(a_1\cdots a_{i-1})|\mathrm{W}(a_{i+1}\cdots a_n)
$$
and the centering property determine then uniquely $\mathrm{W}(a_1\cdots a_n)$. The identities
$$
  \partial_{a_i}\partial_{a_j}\mathrm{W}(a_1\cdots a_n)
  =\partial_{a_j}\partial_{a_i}\mathrm{W}(a_1\cdots a_n)
$$
ensure the consistency of the formulas.
\end{proof}

%%%%%%%%%%%%%%%%%%%%%%%%%%%%%%%%%
%%%%%%%%%%%%%%%%%%%%%%%%%%%%%%%%%

\section{Shuffle algebra}
\label{sect:nsShuffle}

In this section we briefly recall the definition of shuffle algebra, thereby setting the notation used in the rest of the paper. We follow references \cite{EP15,EP18} and refer to these articles for further bibliographical indications on the subject. We use in the present article the topologists' convention and call shuffle products products that are not necessarily commutative (see the definitions below). See also the recent survey \cite{EP19} on the appearance of shuffle algebras (a.k.a.~dendriform algebras) and related structures in the theory of iterated integrals and more generally chronological calculus.

\begin{definition}
A \emph{shuffle algebra} is a vector space $D$ endowed with two bilinear products ${\prec}\colon D\otimes D\to D$ and ${\succ}\colon D\otimes D\to D$, called the left and right half-shuffles, respectively, satisfying the shuffle relations
\begin{equation}
\begin{aligned}
  (a\prec b)\prec c =a\prec(b*c), &
  \quad
  a\succ(b\succ c)=(a*b)\succ c\\
     (a\succ b)\prec c&=a\succ(b\prec c),
  \label{eq:shdef}
\end{aligned}
\end{equation}
where we have set $a*b \coloneq a\succ b+a\prec b$.
\end{definition}
These relations imply that $(D,*)$ is a non-unital associative algebra.
We also consider its unitization $\overline D\coloneq \mathbb{K}\1\oplus D$ by extending the half-shuffles: $\1\prec a\coloneq 0\eqcolon a\succ\1$ and $\1\succ a\coloneq a\eqcolon a\prec\1$ for all $a\in D$. This entails that $1*a=a*1$ for all $a$ in $D$; note however that the products $\1\prec\1$ and $\1\succ\1$ \emph{are not defined}; we put however $\1*\1\coloneq\1$.

\begin{definition}
A commutative shuffle algebra is a shuffle algebra where the left and right half-shuffles are identified by the identity:
\[
  a \succ b - b \prec a=0,
\]
so that in particular $(\overline D,*)$ becomes a commutative algebra and the knowledge of the left-half shuffle $\prec$ (or the right-half shuffle $\succ$) is enough to determine the full structure.
\end{definition}

Shuffle products are frequently denoted $\shuffle$, as we do further below in this article \eqref{eq:comShuffle}. Fundamental examples of such products are provided by the shuffle product of simplices in geometry and topology (see the first part of \cite{EP19} for a modern account) as well as the commutative shuffle product of words defined inductively on $\overline{T}(X)$:
$$
  x_1\cdots x_n\prec y_1\cdots y_m\coloneq x_1(x_2\cdots x_n \shuffle y_1\cdots y_m).
$$
The latter is dual to the unshuffle coproduct $\Delta^\shuffle$. This example is generic in the sense that the tensor algebra over an alphabet $B$ equipped with this product is the free commutative shuffle algebra over $B$ \cite{schu}. The shuffle algebras we will study in the present article are non-commutative variants of the tensor algebra.

Dual to the notion of shuffle algebra is the concept of unshuffle coalgebra \cite{F07}. An \emph{unshuffle coalgebra} is a vector space $C$ equipped with two linear maps $\Delta_\prec\colon C\to C\otimes C$ and $\Delta_\succ\colon C\to C\otimes C$, called the left and right half-unshuffles, such that
\begin{align}
  (\Delta_\prec\otimes{\id})\Delta_\prec&=({\id}\otimes\overline\Delta)\Delta_\prec\label{eq:uca1}\\
  (\Delta_\succ\otimes{\id})\Delta_\prec&=({\id}\otimes\Delta_\prec)\Delta_\succ\label{eq:uca2}\\
  (\overline\Delta\otimes{\id})\Delta_\succ&=({\id}\otimes\Delta_\succ)\Delta_\succ\label{eq:uca3}
\end{align}
where $\overline\Delta \coloneq\Delta_\prec+\Delta_\succ$. As before, these axioms imply that $(C,\overline\Delta)$ is a non-counital coassociative coalgebra.

\begin{definition}\label{def:unshuffle}
  An \emph{unshuffle bialgebra} is a vector space $\overline B=\mathbb{K}\1\oplus B$ together with linear maps $\Delta_\prec\colon B\to B\otimes B$, $\Delta_\succ\colon B\to B\otimes B$ and $m\colon\overline B\otimes\overline B\to\overline B$ such that
\begin{enumerate}
\item $(B,\Delta_\prec,\Delta_\succ)$ is an unshuffle coalgebra,
\item $(\overline B,m)$ is an associative algebra and
\item the following compatibility relations are satisfied:
\end{enumerate}
\[
  \Delta^+_\succ(ab)
  =\Delta^+_\succ(a)\Delta(b),
  \quad
  \Delta^+_\prec(ab)
  =\Delta^+_\prec(a)\Delta(b),
\]
where we have set
\[
  \Delta_\prec^+(a)
  \coloneq\Delta_\prec(a)+a\otimes\1,\quad\Delta_\succ^+(a)\coloneq\Delta_\succ(a)+\1\otimes a
\]
and
\[
  \Delta(a)
  \coloneq\Delta^+_\prec(a)+\Delta^+_\succ(a)
  =\overline\Delta(a)+a\otimes\1+\1\otimes a.
\]
\end{definition}

Given an unshuffle bialgebra, we adjoin a counit $\varepsilon\colon\overline B\to \mathbb{K}$, which is the unique linear map such that $\ker\varepsilon=B$ and $\varepsilon(\1)=1$. We observe that, in particular, for any unshuffle bialgebra the triple $(\overline B,m,\Delta)$ becomes a bialgebra in the usual sense. Thus, its graded dual space $\overline D\coloneq\overline B^*$ becomes an algebra under the convolution product
\begin{equation}
\label{eq:*}
  \varphi * \psi \coloneq(\varphi\otimes\psi)\Delta.
\end{equation}
Moreover, \eqref{eq:uca1}-\eqref{eq:uca3} imply that $\overline D=\mathbb{K}\1 \oplus B^\ast$ is an unital shuffle algebra, since the convolution product splits
\begin{equation}
\label{eq:<>}
  \varphi*\psi=\varphi\prec\psi + \varphi\succ\psi,
\end{equation}
where $\varphi(\1)=\psi(\1)=0$, $\varphi\prec\psi\coloneq(\varphi\otimes\psi)\Delta_\prec^+$ and $\varphi\succ\psi\coloneq(\varphi\otimes\psi)\Delta_\succ^+$. The counit of $\overline B$ plays the role of the unit for this shuffle product and one sets for $\varphi\in \overline D, \ \varphi(\1)=0$,
%\begin{align*}
$$
\begin{array}{c}
  \varepsilon\prec\varphi=(\varepsilon\otimes\varphi)\Delta^+_\prec=0,\\[0.2cm]
  \varphi\succ\varepsilon=(\varphi\otimes\varepsilon)\Delta^+_\succ=0,
\end{array}
\qquad
\begin{array}{c}
  \varphi\prec\varepsilon=(\varphi\otimes\varepsilon)\Delta^+_\prec=\varphi,\\[0.2cm]
  \varepsilon\succ\varphi=(\varepsilon\otimes\varphi)\Delta^+_\succ=\varphi.
\end{array}
$$
By definition, an unshuffle coalgebra is cocommutative if $\tau\circ\Delta_\prec =\Delta_\succ$, where $\tau$ is the usual switch map $\tau(x\otimes y)\coloneq y\otimes x$. An example is given by the algebra $\overline{T}(A)$ equipped with unshuffle coproduct, $\Delta^\shuffle$, defined in \eqref{eq:symcop} below.

\subsection{Shuffle approach to moments and cumulants}
\label{sect:ncrev}
We consider an example of Definition \ref{def:unshuffle}, which is also the main setting for the shuffle algebra approach to moment-cumulant relations in non-commutative probability theory.

We note that the coproduct $\Delta$ can be split into two parts: the \emph{left half-coproduct}
\[
  \Delta_\prec^+(a_1\dotsm a_n)
  \coloneq\sum_{1\in S\subseteq[n]}a_S\otimes a_{J_1^S}\vert\dotsm\vert a_{J_k^S}
\]
and we set
\[
  \Delta_\prec(a_1\dotsm a_n)\coloneq\Delta_\prec^+(a_1\dotsm a_n)-a_1\dotsm a_n\otimes\1.
\]
The \emph{right half-coproduct} is defined by
\begin{equation}
\label{right}
  \Delta_\succ^+(a_1\dotsm a_n)\coloneq\sum_{1\not\in S\subset[n]}a_S\otimes a_{J_1^S}\vert\dotsm\vert a_{J_k^S}
\end{equation}
and we define
\[
  \Delta_\succ(a_1\dotsm a_n)\coloneq\Delta_\succ^+(a_1\dotsm a_n)-\1\otimes a_1\dotsm a_n.
\]
This is extended to the double tensor algebra by defining
\begin{align*}
  \Delta_\prec^+(w_1\vert\dotsm\vert w_m)&\coloneq\Delta_\prec^+(w_1)\Delta(w_2)\dotsm\Delta(w_m)\\
  \Delta_\succ^+(w_1\vert\dotsm\vert w_m)&\coloneq\Delta_\succ^+(w_1)\Delta(w_2)\dotsm\Delta(w_m).
\end{align*}

\begin{theorem}[\cite{EP15}]
  The bialgebra $\overline T(T(A))$ equipped with $\Delta_\prec$ and $\Delta_\succ$ is an unshuffle bialgebra.
\end{theorem}

We recall now from reference \cite{EP15,EP18} how the unshuffle bialgebra $\overline T(T(A))$ provides an algebraic structure for encoding the relation between free, boolean and monotone cumulants and moments in non-commutative probability theory from the point of view of shuffle products.

The group of characters is denoted by $G$ and its Lie algebra of infinitesimal characters $\mathfrak g$ consists of linear maps that send $\1 \in \overline T(T(A))$ as well as any non-trivial product in $ \overline T(T(A))$ to zero. The convolution exponential $\exp^*$ defines a bijection between $\mathfrak g$ and $G$. We recall that the map $\Phi \colon \overline T(T(A))\to \mathbb{K}$ is the unique unital multiplicative extension of the linear map $\phi$ defined on $T(A)$ by $\phi(a_1\dotsm a_n) \coloneq \varphi(a_1\cdot_{\!\scriptscriptstyle{A}} \dotsm\cdot_{\!\scriptscriptstyle{A}} a_n)$. 
We are going to define three different exponential-type bijections between the group $G$ and its Lie algebra $\mathfrak g$, corresponding respectively to the convolution product $*$ and to the right and left half-shuffles (see \Cref{eq:*} and \Cref{eq:<>}). As a result, we can associate to the character $\Phi \in G$ three different infinitesimal characters $\rho, \kappa,\beta \in \mathfrak g$. The three exponential-type bijections encode the three moment-cumulant relations (monotone, free, and boolean). 

The free and boolean cumulants can be represented in terms of infinitesimal characters as the unique maps satisfying the so-called left respectively right half-shuffle fixed point equations
\begin{equation}
  \Phi=\varepsilon+\kappa\prec\Phi
  \qquad
  \mathrm{and}
  \qquad
  \Phi=\varepsilon+\Phi\succ\beta.
  \label{eq:bexp}
\end{equation}
These equations define bijections between the Lie algebra $\mathfrak g$ and the group $G$, i.e., the so-called left and right half-shuffle exponentials such that
\[
    \Phi=\mathcal E_\prec(\kappa)=\mathcal E_\succ(\beta).
\]
Hence, we see that $\Phi$ is the left (or free) half-shuffle exponential of the infinitesimal character $\kappa \in \mathfrak g$. Analogously, $\Phi$ is the right (or boolean) half-shuffle exponential of the infinitesimal character $\beta \in \mathfrak g$. It can be shown \cite[Thm.~5.2]{EP15} that the free moment-cumulant relation of order $n$ is given by computing
$$
  \mathcal E_\prec(\kappa)(a_1 \cdots a_n)=\sum_{\pi \in \operatorname{NC}([n])}\prod_{B\in\pi}\kappa(a_B).
$$
Analogously, $\mathcal E_\succ(\beta)$ gives the boolean moment-cumulant relations \cite[Thm.~4]{EP18} 
$$
  \mathcal E_\succ(\beta)(a_1 \cdots a_n)
  = \sum_{\pi \in \operatorname{Int}([n])} \prod_{B\in\pi} \beta(a_B),
$$
due to the fact that $\beta \in \mathfrak g$ together with the right half-shuffle operation defined in terms of \eqref{right}, implies that
$$
  \mathcal E_\succ(\beta)(a_1 \cdots a_n)
  = \sum_{j=1}^n \Phi(a_{j+1} \cdots a_n) \, \beta(a_1 \cdots a_j).
$$
Shuffle algebra permits to show that half-shuffle exponentials entail the following left and right half-shuffle logarithms
\[
  \kappa=\mathcal{L}_\prec(\Phi)\coloneq(\Phi-\varepsilon)\prec\Phi^{-1},
  \quad
  \beta=\mathcal{L}_\succ(\Phi)\coloneq\Phi^{-1}\succ(\Phi-\varepsilon)
\]
as well as the relation between the boolean and free cumulants through the shuffle adjoint action
\begin{equation}
  \beta=\Theta_\Phi(\kappa)\coloneq\Phi^{-1} \succ\kappa\prec\Phi.
  \label{eq:cumTheta}
\end{equation}

With these notations in place, one can show that the convolutional inverse of $\Phi$ can be also described in terms of the half-shuffle exponentials
\begin{equation}
\label{eq:minvexp}
  \Phi^{-1}=%\exp^*(-\rho)=
  \mathcal E_\succ(-\kappa)=\mathcal E_\prec(-\beta),
\end{equation}
yielding solutions to the half-shuffle fixed point equations
\begin{equation}
  \Phi^{-1}=\varepsilon-\Phi^{-1}\succ \kappa,
  \qquad
  \Phi^{-1}=\varepsilon-\beta\prec\Phi^{-1}.
  \label{eq:bexpinv}
\end{equation}

\subsection{Shuffle calculus for free Wick polynomials}
\label{sect:shW}

The Wick map $\mathrm{W}$ can be related to the free cumulants by using \eqref{eq:minvexp}, whence we obtain
from \Cref{dfn:fWick}
\[
  \mathrm{W}=({\id}\otimes\mathcal E_\succ(-\kappa))\Delta.
\]
Evaluating both sides on a word from $T(A)$ yields
\[
  \mathrm{W}(a_1\dotsm a_n)
  =({\id}\otimes\mathcal E_\succ(-\kappa))\Delta(a_1\dotsm a_n).
\]
Hence from Definition \ref{def:coprod} of the coproduct, we obtain an explicit formula for the Wick polynomial $\mathrm{W}(a_1\dotsm a_n)$, in terms of free cumulants (cf. \cite{Ans04})
\[
  \mathrm{W}(a_1\dotsm a_n)=\sum_{S\subseteq[n]}a_S\sum_{\substack{\pi\in\operatorname{Int}
  ([n]\setminus S)\\\pi\cup S\in \operatorname{NC}([n])}}(-1)^{|\pi|}\prod_{B\in\pi}\kappa(a_B),
\]
which coincides with \cite[Formula (3.44)]{Ans04}.
Note that the combination of the factor $(-1)^{|\pi|}$ and the sum over interval partitions on the right-hand side stems from the fact that $\Phi^{-1}$ is expressed in terms of the right (or boolean) half-shuffle exponential evaluated on the infinitesimal character $-\kappa$ corresponding to negative values of free cumulants. This is the reason for calling these polynomials \emph{free Wick polynomials} and $\mathrm{W}$ is called the \emph{free Wick map}.

\begin{proposition}
  \label{prp:Wickrec}
  The free Wick polynomials satisfy the following recursion in terms of the free cumulants:
  \begin{align}
  \label{freeWrecursion}
  \mathrm{W}=\mathrm{e} + {(\id-\mathrm{e})}\prec\Phi^{-1} - \mathrm{W}\succ\kappa,
  \end{align}
where $\mathrm{e}\coloneq\eta \circ \varepsilon$ and $\eta$ is the unit map on $\overline{T}(T(A))$.
\end{proposition}

\begin{proof}
  This follows from the relations satisfied by the shuffle operations and \eqref{eq:bexpinv}:
  \begin{align*}
    \mathrm{W}&=({\id}\otimes\Phi^{-1})\Delta\\
    &=\mathrm{e}+{(\id-\mathrm{e})}\prec\Phi^{-1}+{\id}\succ(\Phi^{-1}-\varepsilon)\\
    &=\mathrm{e}+{(\id-\mathrm{e})}\prec\Phi^{-1}-{\id}\succ(\Phi^{-1}\succ\kappa)\\
    &=\mathrm{e}+{(\id-\mathrm{e})}\prec\Phi^{-1}-\mathrm{W}\succ\kappa.
  \end{align*}
\end{proof}

We remark that by observing that the left half-coproduct, $\Delta_\prec$, can be expressed in terms of the coproduct $\Delta$, i.e., $\Delta_\prec(a_1\dotsm a_n)=(a_1{\cdot}\otimes{\id})\Delta(a_2\dotsm a_n)$ we recover from \eqref{freeWrecursion} the elegant recursive formula \cite[Formula (3.43)]{Ans04}
\[
  \mathrm{W}(a_1\cdots a_n)=a_1\mathrm{W}(a_2\dotsm a_n)-\sum_{j=0}^{n-1}
  \kappa(a_1\cdots a_j)\,\mathrm{W}(a_{j+1}\dotsm a_n).
\]

\section{Boolean Wick polynomials}
\label{sec:bwickpoly}

It is natural to ask whether one could also relate Wick map, $\mathrm{W}$, to boolean cumulants. Indeed, by using once again \eqref{eq:minvexp}
and \Cref{dfn:fWick}, we obtain
\[
  \mathrm{W}=({\id} \otimes \mathcal E_\prec(-\beta))\Delta.
\]
Expanding the left half-shuffle exponential, $\mathcal E_\prec(-\beta)$, on the righthand side, we see that,
\begin{align*}
  \mathrm{W}
  &={\id}-{\id}\prec\beta-{\id}\succ\beta+{\id}\succ(\beta\prec\beta)+{\id}\prec(\beta\prec\beta)+\dotsb\\
  &=\mathrm{e}+{(\id-\mathrm{e})}\prec(\varepsilon-\beta+\beta\prec\beta+\dotsb)-({\id}\succ\beta)\prec(\varepsilon-\beta+(\beta\prec\beta)+\dotsb)\\
  &=\mathrm{e}+({\id-\mathrm{e}}-{\id}\succ\beta)\prec\Phi^{-1},
\end{align*}
where we have used, in the last identity, the recursion \eqref{eq:bexpinv} and relations \eqref{eq:shdef} to rearrange the iterated half-shuffle products. This argument can be made precise with the help of \Cref{prp:Wickrec}.

\begin{proposition}  \label{prp:bfWick}
  The Wick map can be expressed in terms of boolean cumulants as
  \[
    \mathrm{W}=\mathrm{e}+({\id-\mathrm{e}}-{\id}\succ\beta)\prec\Phi^{-1}.
  \]
\end{proposition}

\begin{proof}
  From \Cref{prp:Wickrec} we have the identity
  \[
  \mathrm{W}=\mathrm{e}+({\id-\mathrm{e}})\prec\Phi^{-1}-\mathrm{W}\succ\kappa
  =\mathrm{e}+({\id-\mathrm{e}})\prec\Phi^{-1}-{\id}\succ(\Phi^{-1}\succ\kappa).
  \]
  But \eqref{eq:cumTheta} implies that $\Phi^{-1}\succ\kappa=\beta\prec\Phi^{-1}$ so that
  \[
    \mathrm{W}=\mathrm{e}+({\id-\mathrm{e}})\prec\Phi^{-1}-{\id}\succ(\beta\prec\Phi^{-1}).
  \]
  Since $a\succ(b\prec c)=(a\succ b)\prec c$ from \eqref{eq:shdef} we get
  \[
    \mathrm{W}=\mathrm{e}+({\id-\mathrm{e}}-{\id}\succ\beta)\prec\Phi^{-1}.
  \]
\end{proof}

We now introduce another map, which allows to recover in a similar way the boolean Appel Polynomials \cite[Section 3]{Ans09a}.
\begin{definition}
  The boolean Wick map $\mathrm{W}'\colon\overline T(T(A))\to\overline T(T(A))$ is defined by
  \begin{equation}
  \label{def:booleanWick}
    \mathrm{W}^\prime\coloneq{\id}-{\id}\succ\beta.
  \end{equation}
\end{definition}

We call as usual boolean Wick polynomials the $\mathrm{W}'(a_1\cdots a_n)$, $a_i\in A$, $i=1,\ldots, n$. In particular, we immediately obtain the explicit expression \cite[Formula (3.1)]{Ans09a}
\begin{equation}
  \mathrm{W}^\prime(a_1\dotsm a_n)
    =a_1\dotsm a_n-\sum_{j=1}^{n}\beta(a_1\dotsm a_{j})a_{j+1}\dotsm a_n.
\label{eq:bwick}
\end{equation}

\begin{proposition}
  The boolean Wick polynomials are centred.
\end{proposition}

\begin{proof}
By definition we have that
\[
    \Phi\circ \mathrm{W}^\prime=\Phi-\Phi\succ\beta=\varepsilon
\]
using \eqref{eq:bexp}.
\end{proof}

\Cref{prp:bfWick} entails the relation
$$
  \mathrm{W}^\prime=\mathrm{e}+(\mathrm{W}-\mathrm{e})\prec\Phi
$$
between the boolean and free Wick maps. This gives the following rewriting rule for the corresponding polynomials,
\begin{equation}
\label{wickprime}
  \mathrm{W}^\prime(a_1\cdots a_n)
  =\sum_{1\in S\subseteq[n]}\mathrm{W}(a_S)\,\Phi(a_{J_1^S})\dotsm\Phi(a_{J_k^S}).
\end{equation}

From \eqref{def:booleanWick} we deduce that
\begin{equation}
\label{idprime}
    {\id}=\mathrm{W}^\prime + {\id}\succ\beta,
\end{equation}
which leads to the expansion
$$
    {\id} = \mathrm{W}^\prime
      + \mathrm{W}^\prime\succ\beta
    + (\mathrm{W}^\prime\succ\beta) \succ \beta
    + ((\mathrm{W}^\prime\succ\beta) \succ\beta) \succ \beta + \cdots.
$$
Observe that the expansion terminates after $n+1$ terms when applied to a word $w\in T(A)$ with $|w|=n$ letters, thanks to $\beta$ being an infinitesimal character, i.e.
\begin{equation}
\label{W'recursion}
  w=\mathrm{W}^\prime(w) 
  	+ \sum_{i=1}^{|w|} R^{(i)}_{\succ\beta}(\mathrm{W}^\prime)(w),
\end{equation}
where $R^{(i)}_{\succ\beta}(\mathrm{W}^\prime)\coloneq R^{(i-1)}_{\succ\beta}(\mathrm{W}^\prime) \succ \beta$ and $R^{(0)}_{\succ\beta}(\mathrm{W}^\prime)=\mathrm{W}^\prime$.
The first few terms are
\begin{align*}
	a_1	   &= \mathrm{W}^\prime(a_1) + \beta(a_1)\\
  a_1a_2      &= \mathrm{W}^\prime(a_1a_2 ) 
  				+ \mathrm{W}^\prime(a_2 ) \beta(a_1) 
				+ \beta(a_1a_2 ) + \beta(a_2) \beta(a_1) ,\\
  a_1a_2a_3     &= \mathrm{W}^\prime(a_1a_2a_3) + \mathrm{W}^\prime(a_2 a_3) \beta(a_1)
          +\mathrm{W}^\prime(a_3 ) \beta(a_1a_2) + \mathrm{W}^\prime(a_3 ) \beta(a_2) \beta(a_1) \\
          	&\quad + \beta(a_1a_2a_3) +  \beta(a_1a_2) \beta(a_3) 
						+ \beta(a_1) \beta(a_2a_3) +  \beta(a_1) \beta(a_2) \beta(a_3)\\
          	&= \mathrm{W}^\prime(a_1a_2a_3) + \mathrm{W}^\prime(a_2 a_3) \beta(a_1)
		+\mathrm{W}^\prime(a_3 ) \big(\beta(a_1a_2) + \beta(a_2) \beta(a_1)\big) 
		+\Phi(a_1 a_2 a_3)   \\	
                &= \mathrm{W}^\prime(a_1a_2a_3) + \sum_{j=1}^{3}\Phi(a_1\cdots a_{j}) \mathrm{W}^\prime(a_{j+1}\cdots a_3)\\
  a_1a_2a_3a_4 &= \mathrm{W}^\prime(a_1a_2a_3a_4)
          + \mathrm{W}^\prime(a_4)\Big(\beta(a_1a_2a_3) + \beta(a_1a_2)\beta(a_3)+ \beta(a_1)\beta(a_2a_3)\\
          &+ \beta(a_1)\beta(a_2)\beta(a_3)\Big)+ \mathrm{W}^\prime(a_3a_4)  \Big(\beta(a_1a_2)
          + \beta(a_1)\beta(a_2)\Big)\\
          &\quad  + \mathrm{W}^\prime(a_2a_3a_4) \beta(a_1) +\Phi(a_1 a_2 a_3a_4).
\end{align*}
Here we used the boolean moment-cumulant relations, which say that $\Phi(a_1 \cdots a_n) = \sum_{I \in\mathrm{Int}([n])} \prod_{\pi \in I}\beta(a_\pi)$.

\begin{proposition}  \label{prp:inverseBoolWick}
	Let $w=a_1\cdots a_n \in T(A)$. Then
\[
    w=\mathrm{W}^\prime(w) +  \sum_{j=1}^{n}\Phi(a_1\dotsm a_{j}) \, \mathrm{W}^\prime(a_{j+1}\dotsm a_n).
\]
\end{proposition}

\begin{proof}
For the word  $w=a_1\dotsm a_n \in T(A)$ we find from \eqref{W'recursion}
\begin{align*}
  w&=\mathrm{W}^\prime(w) + \sum_{j=1}^{n} \mathrm{W}^\prime(a_{j+1} \cdots a_n)
  \sum_{I \in\mathrm{Int}([j])} \prod_{\pi \in I}\beta(a_\pi)\\
    &=\mathrm{W}^\prime(w) +  \sum_{j=1}^{n} \Phi(a_{1}\cdots a_j)\,\mathrm{W}^\prime(a_{j+1} \cdots a_n),
	\end{align*}
The essential input here is that the boolean cumulants are given by $\beta$, which is an infinitesimal character.
\end{proof}

Eventually, from \eqref{idprime} we deduce the inverse boolean Wick map.

\begin{proposition}  \label{prp:inversebooleanWick} 
The inverse boolean Wick map is given as solution to the fixed point equation
\begin{equation}
\label{W'inverse}
	\mathrm{W}'^{\circ -1} = \id + \mathrm{W}'^{\circ -1} \succ \beta.
\end{equation}
\end{proposition}

\begin{proof}
Note that the definition of the boolean Wick map \eqref{def:booleanWick} implies that it is invertible. We show explicitly that  
$\mathrm{W}'^{\circ -1} \circ \mathrm{W}' = \mathrm{W}' \circ \mathrm{W}'^{\circ -1} = \id.$ Indeed, we see that 
$$
	 \mathrm{W}' \circ \mathrm{W}'^{\circ -1} =  \mathrm{W}' +  (\mathrm{W}'\circ  \mathrm{W}'^{\circ -1}) \succ \beta.
$$
Induction on the length of words in $T(A)$ gives for $a \in A$
$$
	 \mathrm{W}' \circ \mathrm{W}'^{\circ -1}(a) =  \mathrm{W}'(a) +  \beta(a) = a.
$$
On a word $w=a_1\dotsm a_n \in T(A)$, $n>1$, we find
\begin{align*}
	 \mathrm{W}' \circ \mathrm{W}'^{\circ -1}(w) 
	 &=  \mathrm{W}'(w) 
	 	+  \sum_{i=1}^{n}(\mathrm{W}' \circ \mathrm{W}'^{\circ -1}) (a_{i+1} \cdots a_n) \beta(a_{1} \cdots a_{i}) \\
	 &=  \mathrm{W}'(w) +  \sum_{i=1}^{n}a_{i+1} \cdots a_n \beta(a_{i+1} \cdots a_n) \\
	 &=w
\end{align*}
Here we used the induction hypothesis, $(\mathrm{W}' \circ \mathrm{W}'^{\circ -1}) (a_{i+1} \cdots a_n) = a_{i+1} \cdots a_n$, for  $i>0$. An analogue computation gives the opposite, i.e., $\mathrm{W}'^{\circ -1} \circ \mathrm{W}' =\id$.
\end{proof}

From \eqref{W'inverse} it follows that
$$
	\mathrm{W}'^{\circ -1}(a_1\dotsm a_n )
	= \sum_{j=0}^n \Phi(a_1\dotsm a_j)a_{j+1}\dotsm a_n  
$$

\begin{remark}
In \cite{Ans09a} the boolean cumulants were defined by the relation between generating functions $G$ of boolean Wick polynomials and boolean cumulants $\eta$,
\[
    G(x,z)=(1-x\cdot z)^{-1}(1-\eta(z)),
\]
which implies an expression similar to \eqref{eq:bwick} but with $\beta$ applied to the other half of the word. In principle, one could take either relation as a starting point since there is a choice here due to the non-commutativity of the series, and neither choice seems to be more natural than the other. However, we decided to work with \eqref{eq:bwick} instead since the polynomials so obtained are more naturally described from the shuffle algebra point of view. The relation (\ref{wickprime}) also has its counterpart in terms of generating functions, which involves a particular kind of variable substitution.
\end{remark}

%%%%%%%%%%%%%%%%%%%%%%%%%%%%%%%%%
%%%%%%%%%%%%%%%%%%%%%%%%%%%%%%%%%

\section{Conditionally free Wick polynomials}
\label{sect:cfreeW}

Note the apparent asymmetry in the definitions of the free and boolean Wick polynomials. There is a third family of polynomials that generalises both the free and boolean cases. Indeed, we may consider the notion of conditional freeness \cite{Boz1996} which generalises Voiculescu's notion of freeness in the context of two states. Recall that a \emph{two-state non-commutative probability space} \((A,\varphi,\psi)\) is a non-commutative probability space $(A,\varphi)$ endowed with a second unital linear map $\psi\colon A \to \mathbb{K}$. We denote by $\Psi$ the canonical character extension of $\psi$ to the double tensor algebra $\overline T(T(A))$. We denote by $\beta^\varphi$ the boolean infinitesimal character associated to $\varphi$ (and define similarly $\beta^\psi,\kappa^\varphi,\kappa^\psi$).

In the shuffle algebra approach we have the following characterisation of conditionally (c-)free cumulants \cite{EP18a}: the corresponding infinitesimal character $R^{\varphi,\psi} \in \mathfrak{g}$ is defined through shuffle adjoint action:
\begin{equation}
  \label{eq:cfreebool}
  R^{\varphi,\psi}\coloneq\Psi\succ\beta^\varphi\prec\Psi^{-1}.
\end{equation}
This means that $\beta^\varphi=\Psi^{-1}\succ R^{\varphi,\psi}\prec\Psi$, such that
\begin{equation}
  \label{eq:cfreeboolGOOD}
    \Phi=\varepsilon+\Phi\succ(\Psi^{-1}\succ R^{\varphi,\psi}\prec\Psi).
\end{equation}

Following \cite[Prop.~6.1]{EP18a} the evaluation of formula \eqref{eq:cfreeboolGOOD} on a word, i.e., computing $\Phi(a_1 \cdots a_n)
	= \varphi(a_1\cdot_{\!\scriptscriptstyle{A}} \dotsm\cdot_{\!\scriptscriptstyle{A}} a_n)$ 
$$
	\Phi(a_1 \cdots a_n)
	=\Phi\succ(\Psi^{-1}\succ R^{\varphi,\psi}\prec\Psi)(a_1 \cdots a_n)
$$ 
gives back the formula discovered in reference \cite{Boz1996} and recalled in the next theorem.

\begin{theorem}[\cite{Boz1996}]
  The following relation between moments and conditionally free cumulants holds:
\begin{equation}
  \label{eq:cfreecumu}
    \varphi(a_1\cdot_{\!\scriptscriptstyle{A}} \dotsm\cdot_{\!\scriptscriptstyle{A}} a_n)
    =\sum_{\pi\in \operatorname{NC}([n])}\prod_{B\in\operatorname{Outer}(\pi)}
  R^{\varphi,\psi}(a_B)\prod_{B\in\operatorname{Inner}(\pi)}\kappa^\psi(a_B).
\end{equation}
\end{theorem}
Here, a block $\pi_i$ of a non-crossing partition $\pi \in NC_n$ is ``Inner'' if there exists a $\pi_j$ and $a,b\in \pi_j$ such that $a<c<b$ for all $c\in \pi_i$. A block which is not an inner one is ``Outer''.

Conditionally free cumulants contain both free and boolean cumulants as limiting cases. More precisely, if we consider the case $\psi=\varphi$ then \eqref{eq:cfreebool} entails
\[
  R^{\varphi,\varphi}
  =\Phi\succ\beta^\varphi\prec\Phi^{-1}
  =\kappa^\varphi
\]
by \eqref{eq:cumTheta}. On the other hand, if $\psi=\varepsilon$ is the trivial state then
\[
  R^{\varphi,\varepsilon}=\beta^\varphi.
\]

\begin{theorem}[\cite{EP17}]
  Let $\alpha_1,\alpha_2$ be two infinitesimal characters of the double tensor algebra and denote by $\mathcal E_\succ(\alpha_1)$ and $\mathcal E_\succ(\alpha_2)$ the corresponding right half-shuffle exponentials. The right half-shuffle Baker--Campbell--Hausdorff formula holds:
  \[
    \mathcal{L}_\succ(\mathcal E_\succ(\alpha_1)*\mathcal E_\succ(\alpha_2))
    =\alpha_2+\Theta_{\mathcal E_\succ(\alpha_2)}(\alpha_1),
   \]
   \label{thm:shlog}
   where $\Theta$ stands for the (shuffle) adjoint action:
$$
  \Theta_{\mathcal E_\succ(\alpha_2)}(\alpha_1)
  = \mathcal E^{-1}_\succ(\alpha_2)\succ \alpha_1\prec \mathcal E_\succ(\alpha_2).
$$
\end{theorem}

\begin{proof}
  Let $X=\mathcal E_\succ(\alpha_1)$ and $Y=\mathcal E_\succ(\alpha_2)$.
  By definition of the shuffle product we have that
  \begin{align*}
    X*Y-\varepsilon&=(X-\varepsilon)\prec Y+X\succ (Y-\varepsilon)\\
    &=(X\succ \alpha_1)\prec Y+X\succ(Y\succ \alpha_2)\\
    &=(X\succ \alpha_1)\prec Y+(X*Y)\succ \alpha_2.
  \end{align*}
  Now, observe that
  \begin{align*}
    (X\succ \alpha_1)\prec Y&=((X*Y*Y^{-1})\succ \alpha_1)\prec Y\\
    &=(X*Y\succ (Y^{-1}\succ \alpha_1))\prec Y\\
    &=X*Y\succ (Y^{-1}\succ \alpha_1 \prec Y).
  \end{align*}
This implies the result using the definition of $\mathcal{L}_\succ$.
\end{proof}

Returning to \Cref{dfn:fWick}, since $\Phi=\mathcal E_\succ(\Theta_{\Psi}(R^{\varphi,\psi}))$ and $\Phi^{-1}=\mathcal E_\prec(-\Theta_{\Psi}(R^{\varphi,\psi}))$
we may now express the free Wick map $\mathrm{W}=\mathrm{id}\ast \Phi^{-1}$ in terms of the conditionally free cumulants $R^{\varphi,\psi}$ as
\[
  \mathrm{W}=\big({\id} \otimes \mathcal E_\prec(-\Theta_{\Psi}(R^{\varphi,\psi}))\big)\Delta.
\]
A computation similar to the boolean case yields
\begin{align*}
  \mathrm{W}&=\big({\id} \otimes \mathcal E_\prec(-\Theta_{\Psi}(R^{\varphi,\psi}))\big)\Delta\\
     &=\mathrm{e}+\big(\id - \mathrm{e} - \id \succ \Theta_{\Psi}(R^{\varphi,\psi}) \big)\prec \Phi^{-1}\\
     &=\mathrm{e}+\big((\id - \mathrm{e}) \prec \Psi^{-1} - \id \succ( \Psi^{-1}  \succ R^{\varphi,\psi} \big)\big)\prec (\Phi*\Psi^{-1})^{-1}.
\end{align*}

\begin{definition}
The conditionally free Wick polynomials are defined to be
\begin{equation}
\label{eq:cfreeWick}
  \mathrm{W}^c  \coloneq \mathrm{e}+(\mathrm{W} - \mathrm{e}) \prec \Phi*\Psi^{-1},
\end{equation}
\end{definition}

This means
\begin{align*}
  \mathrm{W}^c
  &\coloneq \mathrm{e}+(\id-\mathrm{e}) \prec \Psi^{-1} - \id \succ( \Psi^{-1}  \succ R^{\varphi,\psi})\\
      &= \mathrm{e}+(\id-\mathrm{e})  \prec \Psi^{-1} - (\id * \Psi^{-1})  \succ R^{\varphi,\psi}\\
      &= \mathrm{e} +\big(\id - \mathrm{e} - \id \succ \Theta_{\Psi}(R^{\varphi,\psi})\big)\prec \Psi^{-1}
\end{align*}

From \eqref{eq:cfreeWick} we deduce a --intricate-- recursion for the inverse of conditionally free Wick map:
\begin{align}
\label{invWcond}
\begin{split}
	\lefteqn{{\mathrm{W}^c}^{\circ -1}(a_1 \cdots a_n)
	= a_1 \cdots a_n - {\mathrm{W}^c}^{\circ -1}\circ(\mathrm{W} - \id)(a_1 \cdots a_n) }	\\
	&\qquad
	 - \sum_{1\in S \subsetneq [n]} {\mathrm{W}^c}^{\circ -1}\circ\mathrm{W}(a_S)
	 (\Phi*\Psi^{-1})(a_{J^S_1}) \dotsm (\Phi*\Psi^{-1})(a_{J_l^S})
\end{split}
\end{align}

Starting again from the identity
\[
  \mathrm{W}={\id}*\Phi^{-1}={\id}*\mathcal E_\prec(-\Theta_\Psi(R^{\varphi,\psi})),
\]
 we obtain after some simple manipulations
\begin{align*}
  \mathrm{W}
  &={\id}*\Psi^{-1}*\Psi*\mathcal E_\prec(-\Theta_\Psi(R^{\varphi,\psi}))\\
  &=\mathrm{W}^\psi*\mathcal E_\prec(\kappa^\psi)*\mathcal E_\prec(-\Theta_\Psi(R^{\varphi,\psi}))\\
  &=\mathrm{W}^\psi*\mathcal E_\prec(\kappa^\psi-R^{\varphi,\psi}),
\end{align*}
where we have used \Cref{thm:shlog} in the last equality. Hence, we have that the free Wick maps $\mathrm{W}$ and $\mathrm{W}^\psi\coloneq(\id * \Psi^{-1})  $ are related
\begin{align*}
  \mathrm{W}
    &= \mathrm{W}^\psi * (\Psi* \Phi^{-1})\\
        &= \mathrm{W}^\psi*\mathcal E_\prec(\kappa^\psi-R^{\varphi,\psi}).
\end{align*}

Finally, we observe from \eqref{eq:cfreeWick} that, in the cases $\Psi=\Phi$ and $\Psi=\varepsilon$ we recover the free and boolean Wick maps $\mathrm{W}$ and $\mathrm{W}^\prime$, respectively.

%%%%%%%%%%%%%%%%%%%%%%%%%%%%%%%%%
%%%%%%%%%%%%%%%%%%%%%%%%%%%%%%%%%

\section{Wick polynomials as group actions}
\label{sect:group}

Observe that the coproduct defined in \cref{def:coprod} is linear on the left and polynomial on the right factor when restricted to $T(A)$, i.e., $\Delta\colon T(A)\to \overline T(A)\otimes\overline T(T(A))$. This means in particular that $\overline T(A)$ is a right comodule over $\overline T(T(A))$, simply by coassociativity. Thus we can induce an action of the group $G$ of characters over $\overline T(T(A))$ on the space $\operatorname{End}(\overline T(A))$ of linear endomorphisms of $\overline T(A)$ by setting
\[
  L . \Psi=(L\otimes \Psi)\Delta.
\]
More precisely we have

\begin{proposition}
  Given $\Psi \in G$ and $L \in\operatorname{End}(\overline T(A))$, define $L . \Psi \in\operatorname{End}(\overline T(A))$ as above.
Then $(\Psi ,L) \mapsto L.\Psi$ defines a (right) action of $G$ on $\operatorname{End}(\overline T(A))$.
\end{proposition}

\begin{proof}
  Let $\Psi_1,\Psi_2 \in G$ and $L \in \operatorname{End}(\overline T(A))$.
  Clearly $L.\Psi \in \operatorname{End}(\overline T(A))$ and
  \begin{align*}
    (L.\Psi_1).\Psi_2
    &= (L\otimes \Psi_1 \otimes \Psi_2) \circ(\Delta \otimes {\id})\circ\Delta\\
    &= (L\otimes \Psi_1\otimes \Psi_2)\circ({\id}\otimes\Delta)\circ\Delta\\
    &= (L\otimes \Psi_1*\Psi_2)\circ\Delta\\
    &= L.(\Psi_1*\Psi_2)
  \end{align*}
  so the mapping $(\Psi_1,L)\mapsto L.\Psi_1$ is an action of $G$ on $\operatorname{End}(\overline T(A))$.
\end{proof}

In the following we identify implicitly the (various) notions of Wick polynomials with the (various) restrictions of the Wick maps to $\overline T(A)$. So, in this section and the following, $\mathrm{W}$ denotes the restriction of $\mathrm{W}$ to $\overline T(A)$, and so on (as should be anyway clear from the context).

As we have seen above, the orbit of the identity map ${\id}\in\operatorname{Aut}(\overline T(A))$ consists only of automorphisms of $\overline T(A)$ and we have the inversion formula for the composition of endomorphisms $({\id}.\Psi)^{-1}={\id}.\Psi^{-1}$ where on the right hand side $\Psi$ is inverted with respect to convolution.
The free Wick polynomials $\mathrm{W}=\id.\Phi^{-1}$ are elements in the orbit of the identity endomorphism by the group action of $G$ on $\operatorname{End}(\overline T(A))$.

Regarding the left half-unshuffle coproduct $\Delta_\prec^+$, we get from \eqref{eq:uca1} that $(T(A),\Delta_\prec)$ is also a right-comodule over $(\overline T(T(A)), \Delta)$. At the level of endomorphisms, we obtain the following

\begin{proposition}
  Let $L\in\operatorname{End}(T(A))$ and $\Psi \in G$. The composition $(\Psi,L)\mapsto L^\Psi \coloneq(L\otimes \Psi)\Delta_\prec$ defines a (right) action.
\end{proposition}

Thus we might reinterpret the boolean Wick polynomials $\mathrm{W}^\prime=\mathrm{e}+(\mathrm{W}-\mathrm{e})\prec \Phi$ as being given on $T(A)$ by a combined action $\mathrm{W}^\prime=\mathrm{e}+({\id}.\Phi^{-1}-\mathrm{e})^{\Phi}$. More generally the relation between the conditionally free and free Wick polynomials can be re-expressed on $T(A)$ as
\[
  \mathrm{W}^c
  =\mathrm{e}+({\id}.\Phi^{-1}-\mathrm{e})^{\Phi*\Psi^{-1}}
  =\mathrm{e}+\left[ \left( {\id}.\Phi^{-1} -\mathrm{e}\right)^{\Phi} \right]^{\Psi^{-1}}.
\]

Neglecting the degree zero (that is, the $\mathrm{e}$) terms, the relations between free, boolean and conditionally free Wick polynomials are encoded by the following diagram:
\begin{figure}[!ht]
  \centering
  \begin{tikzcd}
    {\id}\rar[".\Phi^{-1}"]&\mathrm W \rar["()^{\Phi}"] &\mathrm W'\rar["()^{\Psi^{-1}}"] &\mathrm W^c.
  \end{tikzcd}
  \label{fig:wickact}
\end{figure}

%%%%%%%%%%%%%%%%%%%%%%%%%%%%%%%%%
%%%%%%%%%%%%%%%%%%%%%%%%%%%%%%%%%

\section{Free, boolean and conditionally free Wick products}
\label{sect:products}

Let $(A,\varphi)$ be a non-commutative probability space. Let $F \colon \overline T(A)\to \overline T(A)$ be an invertible linear map such that $F(1_A)=1_A$. One can induce a modified product ${\bullet}$ on $\overline T(A)$ by conjugacy, that is setting $w\bullet w' \coloneq F(F^{-1}(w)F^{-1}(w'))$. Associativity follows from associativity of the concatenation product on $\overline T(A)$. Therefore, $F$ becomes a unital algebra morphism from $(\overline T(A),\otimes)$ to $(T(A),\bullet)$.

Since the maps $\mathrm{W}$, $\mathrm{W}^\prime$ and $\mathrm{W}^c$ are all invertible when acting on $\overline T(A)$ we obtain from this construction three new products on $\overline T(A)$:

\begin{definition}
The three associative products on $\overline T(A)$ induced by the three Wick maps $\mathrm{W}$, $\mathrm{W}^\prime$ and $\mathrm{W}^c$ are denoted by $\bullet$, $\odot$ and $\times$ and called the free, boolean and conditionally free Wick products, respectively. The Wick maps are morphisms of algebras when $\overline T(A)$ is equipped with either of these new products. In particular for $a\in A$,
$$
  \mathrm{W}(a^n)=\mathrm{W}(a)^{\bullet n},
$$
and similarly for the other cases.
\end{definition}

The conjugacy formula gives the rule for computing the new products. For example, in the free and Boolean case we find

\begin{proposition}
\label{prp:deformedproducts}
\begin{enumerate}
\item  The free Wick product $\bullet$ admits the following closed-form formula: for words $w=a_1\dotsm a_n$ and $w'=a_{n+1}\dotsm a_{n+m}$ in $T(A)$, we find
\[
    w\bullet w'=\sum_{S\subseteq[n+m]}\mathrm{W}(a_S)\Phi(a_{K^S_1})\dotsm\Phi(a_{K_l^S}).
\]
where the $K_i^S,i=1,\dots , l$ run over the connected components of $[n]-([n]\cap S)$ and $(n+[m])-(n+[m]\cap S)$.
   
\item The boolean Wick product $\odot$ admits the following closed-form formula: for words $w=a_1\dotsm a_n$ and $w'=b_{1}\dotsm b_{m}$ in $T(A)$, we find 
\[
		w \odot w'= \sum_{\substack{0 \leq i \leq n\\0 \leq  j \leq m}
    	\Phi(a_1 \cdots a_i \vert b_{1} \cdots b_j) \mathrm{W}'(a_{i+1} \cdots a_n b_{j+1} \cdots b_m).}
\] 
\end{enumerate}
\end{proposition}

\begin{proof} 
\begin{enumerate}
\item
Set $b_i\coloneq a_{n+i}, \ i=1,\dots , m$.
  Since the inverse free Wick map is the map $\mathrm{W}^{\circ -1}=({\id}\otimes\Phi)\Delta$ we have that
  \begin{align*}
    \mathrm{W}^{\circ -1}(w)\mathrm{W}^{\circ -1}(w')
    &=\sum_{S\subseteq[n]}\sum_{S'\subseteq [m]}
    a_S\,b_{S'}\,\Phi(a_{J_1^S})\dotsm\Phi(a_{J_{k(S)}^S})\Phi(b_{J_1^{S'}})\dotsm\Phi(b_{J_{k(S')}^{S'}}).
  \end{align*}
By re-expressing in terms of the $a_i$ we get
$$
  \mathrm{W}^{\circ -1}(w)\mathrm{W}^{\circ -1}(w')
  =\sum_{S\subseteq[n+m]}a_S\,\Phi(a_{K^S_1})\dotsm\Phi(a_{K_l^S}).
$$
The conclusion then follows by applying $\mathrm{W}$ to both sides of this identity.

\item
 Recall Proposition \ref{prp:inversebooleanWick} stating that the inverse boolean Wick map is given recursively $\mathrm{W}'^{\circ -1}=\id + \mathrm{W}'^{\circ -1} \succ \beta$ such that
$$
	 \mathrm{W}'^{\circ -1}(a_1\dotsm a_n)
	=  \sum_{j=0}^n \Phi(a_1\dotsm a_j) a_{j+1}\dotsm a_n.
$$ 
Then we have
$$
	\mathrm{W}'^{\circ -1}(a_1\dotsm a_n)
	\mathrm{W}'^{\circ -1}(b_1\dotsm b_m)=
	 \sum_{\substack{0 \leq i \leq n\\
	 0 \leq  j \leq m}}
    	\Phi(a_1 \cdots a_i \vert b_{1} \cdots b_j) a_{i+1} \cdots a_n b_{j+1} \cdots b_m.
$$
The conclusion then follows by applying $\mathrm{W}'$ to both sides of this identity.	
\end{enumerate}
\end{proof}

\begin{remark}
A closed formula for the conditionally free Wick products follows from using the recursion \eqref{invWcond}. 
\begin{align*}
\lefteqn{{\mathrm{W}^c}^{\circ -1}(a_1 \cdots a_n){\mathrm{W}^c}^{\circ -1}(b_1 \cdots b_m)}\\
	&= \Big(a_1 \cdots a_n - {\mathrm{W}^c}^{\circ -1}\circ(\mathrm{W} - \id)(a_1 \cdots a_n)\\
	&\qquad
	 - \sum_{1\in S \subsetneq [n]} {\mathrm{W}^c}^{\circ -1}\circ\mathrm{W}(a_S)
	 (\Phi*\Psi^{-1})(a_{J^S_1}) \dotsm (\Phi*\Psi^{-1})(a_{J_l^S})\Big)\\
	& \Big(b_1 \cdots b_m - {\mathrm{W}^c}^{\circ -1}\circ(\mathrm{W} - \id)(b_1 \cdots b_m)\\	
	&\qquad
	 - \sum_{1\in S \subsetneq [m]} {\mathrm{W}^c}^{\circ -1}\circ\mathrm{W}(b_S)
	 (\Phi*\Psi^{-1})(b_{J^S_1}) \dotsm (\Phi*\Psi^{-1})(b_{J_l^S})\Big)
\end{align*}
Applying ${\mathrm{W}^c}$ on both sides gives the conditionally free Wick product 
$$
	a_1 \cdots a_n \times b_1 \cdots b_m
	={\mathrm{W}^c}\big({\mathrm{W}^c}^{\circ -1}(a_1 \cdots a_n){\mathrm{W}^c}^{\circ -1}(b_1 \cdots b_m)\big).
$$
%However, for the sake of space limit we refrain from stating it explicitly.}
\end{remark}

%%%%%%%%%%%%%%%%%%%%%%%%%%%%%%%%%
%%%%%%%%%%%%%%%%%%%%%%%%%%%%%%%%%

\section{Tensor cumulants}
\label{ssect:tensorc}
We now briefly show how our approach allows to lift the classical notion of cumulants to the non-commutative setting and to revisit the notion of tensor cumulants \cite{NS06} as a warm up for the definition of tensor Wick polynomials.

As before we work on a non-commutative probability space \( (A,\varphi)\) (see \cref{dfn:ncproba}).
On $\overline T(A)$  the unshuffle coproduct $\Delta^{\shuffle}\colon \overline T(A)\to \overline T(A)\otimes \overline T(A)$ is defined by declaring elements in $A \hookrightarrow \overline T(A)$ to be primitive and extending it multiplicatively to all of $\overline T(A)$. As a result one gets that for any $a_1,\dotsc, a_n \in A$,
\begin{align}
  \Delta^\shuffle(a_1\dotsm a_n)
      &=\Delta^\shuffle(a_1)\dotsm \Delta^\shuffle(a_n)     \nonumber \\
                &= \sum_{S\subseteq[n]}a_S\otimes a_{[n]\setminus S},
  \label{eq:symcop}
\end{align}
where we have set $a_\emptyset\coloneq\1$ and
$$
  a_U\coloneq{a_{u_1}\dotsm a_{u_p}}
$$
for $U=\{u_1 < \cdots < u_p\}\subseteq[n]$. This endows the unital tensor algebra with the structure of a cocommutative graded connected Hopf algebra. The antipode reverses the order of the letters in a word and multiplies it by a minus sign if the word has odd length.

Its dual  $\overline T(A)^*$ is a commutative algebra with the convolution product defined for linear maps $\mu,\nu \colon \overline T(A) \to \mathbb{K}$ by the commutative shuffle product
\begin{equation}
\label{eq:comShuffle}
  \mu\shuffle\nu \coloneq (\mu\otimes\nu)\Delta^{\shuffle}.
\end{equation}
The unit for this product is the counit $\varepsilon\colon\overline T(A) \to \mathbb{K}$, which is uniquely defined by $\ker\varepsilon=T(A)$ and $\varepsilon(\1)=1$. See e.g. \cite{reutenauer1993free} for details.

The generalized expectation map $\varphi$ permits to define a linear map $\phi \colon \overline T(A)\to \mathbb{K}$ by setting $\phi(a_1\dotsm a_n) = \varphi(a_1\cdot_{\!\scriptscriptstyle{A}} \dotsm\cdot_{\!\scriptscriptstyle{A}} a_n)$ and $\phi(\1)=1$.

The grading on $\overline T(A)$ permits to think of $\phi$ as a graded series
\[
  \phi=\sum_{n\ge 0}\phi_n,
\]
where $\phi_n\colon T(A) \to \mathbb{K}$ is a linear map vanishing outside $T_n(A)$, the degree $n$ component of $T(A)$. In this way, we may regard the map $\phi$ as being some kind of generalized moment-generating function. Since the algebra $T(A)$ is graded by the length of words and connected ($T_0(A)=\mathbb K\1$), the exponential and logarithm maps define inverse bijections between unital linear maps on $\overline{T}(A)$ and reduced maps (maps that vanish on $\mathbb{K}$, the degree zero component). In particular  there exists a unique linear map $c \in \overline T(A)^*$ with $c(\1)=0$ such that
$$
  \phi=\exp^\shuffle(c),\quad\ \phi^{-1}=\exp^\shuffle(-c),
$$
where $\phi^{-1}$ is the inverse of $\phi$ for the shuffle product \eqref{eq:comShuffle}.

\begin{definition}
The tensor cumulant map associated to $\phi $ is the linear application $c\colon\overline T(A) \to \mathbb{K}$ defined by
$$
  c\coloneq\log^\shuffle(\phi).
$$
Its evaluations $c(a_1\cdots a_n) \in \mathbb{K}$ are also written $c(a_1,\dots ,a_n)$ and are called the multivariate tensor cumulants associated to the sequence $(a_1,\dots,a_n)$ of non-commutative random variables.
\end{definition}

The defining relation $c\coloneq\log^\shuffle(\phi)$ is a version in a non-commutative context of the usual formula relating the moment and cumulant generating functions (see \eqref{eq:cmcrel}). From \eqref{eq:symcop} we see that for any $j>0$ the iterated reduced coproduct $\overline\Delta^\shuffle_{j-1}\colon T(A) \to T(A)^{\otimes j}$ is given by $\overline\Delta^\shuffle_{0}=\id$ and
\begin{equation}
  \overline\Delta^\shuffle_{j-1}(a_1\dotsm a_n)
      =\sum_{\pi\in P_j(n)}\sum_{\sigma\in\mathbb S_j}a_{B_{\sigma(1)}}\otimes\dotsm\otimes a_{B_{\sigma(j)}},
  \label{eq:copiter}
\end{equation}
where $P_j(n)$ is the collection of all set partitions $\pi=\{B_1,\ldots,B_j\}$ of $[n]\coloneq\{1,\dotsc,n\}$ into $j$ disjoint subsets and $\mathbb S_j$ is the $j$-th symmetric group (recall that for $x\in T(A)$, $\overline\Delta^\shuffle(x)\coloneq\Delta^\shuffle(x)-x\otimes 1-1\otimes x$). From $c(\1)=0$ we deduce
$$
  \phi=\exp^\shuffle(c)=c
        +\frac12(c\otimes c)\overline\Delta_1^\shuffle
        +\frac16(c\otimes c\otimes c)\overline\Delta^\shuffle_2
        +\frac{1}{24}(c\otimes c\otimes c\otimes c)\overline\Delta^\shuffle_3
        +\dotsb,
$$
giving the multidimensional version of formula \eqref{eq:cmcrel}
\begin{equation}
  \phi(a_1 \dotsm a_n)
  = \varphi(a_1\cdot_{\!\scriptscriptstyle{A}} \dotsm\cdot_{\!\scriptscriptstyle{A}} a_n)
  = \sum_{\pi\in P(n)}\prod_{B\in\pi}c (a_B).
  \label{eq:mcrel}
\end{equation}
Recall that $c (a_B):=c (a_{i_1}, \ldots, a_{i_{|B|}})$, for $B=\{i_1 < \cdots < i_{|B|}\}$, is the multivariate cumulant of order $|B|$ and $P(n)$ is the collection of all set partitions of $[n]$. In fact, many other versions of this relation can be recovered from the properties of the underlying Hopf algebra. See \cite{EP15,EPTZ18} for further details. The important point here is that set partitions appear naturally in \eqref{eq:mcrel} through formula \eqref{eq:copiter} due to the definition of the coproduct in \eqref{eq:symcop}.

%%%%%%%%%%%%%%%%%%%%%%%%%%%%%%%%%
%%%%%%%%%%%%%%%%%%%%%%%%%%%%%%%%%

\subsection{Tensor Wick polynomials}
\label{sect:Wick}

It turns out that the same Hopf-algebraic framework used for describing the tensor moment-cumulant relations allows to get an explicit description of tensor Wick polynomials (that can be understood as a natural non-commutative lift of classical Wick polynomials).

\begin{definition} The \emph{tensor Wick map} $W_T\colon \overline T(A)\to \overline T(A)$ is defined by
\[
  W_T\coloneq ({\id}\otimes\phi^{-1})\Delta^\shuffle.
\]
Its inverse $W_T^{-1}\colon \overline T(A)\to \overline T(A)$ is given by
\[
  W_T^{-1} = ({\id}\otimes\phi)\Delta^\shuffle.
\]
Given a sequence $(a_1,\dots,a_n)\in A$, $W_T(a_1\cdots a_n)$ is called the tensor Wick polynomial associated to this sequence.
\end{definition}

Let us compute a few examples using the reduced unshuffle coproduct \eqref{eq:copiter} and the fact that the inverse $\phi^{-1} \colon \overline T(A)\to \mathbb{K}$ is given by the Neumann series
\[
  \phi^{-1}=\varepsilon + \sum_{n>0}(-1)^n \phi^{\otimes n}\overline\Delta^\shuffle_{n-1}.
\]
Then the first three tensor Wick polynomials in $\overline T(A)$ are
\begin{align}
  W_T(a_1)&=a_1 - \varphi(a_1)\1 \nonumber\\
  W_T(a_1a_2)&=a_1a_2 - a_2\varphi(a_1) - a_1\varphi(a_2)
  + \big(-\varphi(a_1 \cdot_{\!\scriptscriptstyle{A}} a_2)
    + 2 \varphi(a_1)\varphi(a_2)\big)\1 \nonumber \\
  \begin{split}
  W_T(a_1a_2a_3)&=a_1a_1a_3 - a_2a_3\varphi(a_1) - a_1a_3\varphi(a_2) - a_1a_2\varphi(a_3)
   +a_1\big(-\varphi(a_2 \cdot_{\!\scriptscriptstyle{A}} a_3) \\
  & + 2 \varphi(a_2)\varphi(a_3)\big) + a_2(-\varphi(a_1 \cdot_{\!\scriptscriptstyle{A}} a_3) + 2 \varphi(a_1)\varphi(a_3))
      +a_3 \big(-\varphi(a_1 \cdot_{\!\scriptscriptstyle{A}} a_2)\\
  &      + 2 \varphi(a_1)\varphi(a_2)\big) + \big(- \varphi(a_1 \cdot_{\!\scriptscriptstyle{A}} a_2 \cdot_{\!\scriptscriptstyle{A}} a_3)
      + \varphi(a_1)\varphi(a_2 \cdot_{\!\scriptscriptstyle{A}} a_3)\\
  & + \varphi(a_2)\varphi(a_1\cdot_{\!\scriptscriptstyle{A}} a_3)
    + \varphi(a_3)\varphi(a_1\cdot_{\!\scriptscriptstyle{A}} a_2)
      - 6 \varphi(a_1)\varphi(a_2)\varphi(a_3)\big)\1.
  \label{Wick3}\end{split}
\end{align}

\begin{remark}\label{rmk:properWick} The tensor Wick map $W_T$ associates to words  $w \in \overline T(A)$ non-commutative polynomials $W_T(w)$ in $\overline T(A)$. Saying this, if the algebra $A$ is commutative then those non-commutative polynomials map by the evaluation $ev\colon\ a_1 \cdots a_n\longmapsto a_1 \cdot_{\!\scriptscriptstyle{A}} \dots \cdot_{\!\scriptscriptstyle{A}} a_n$ to the  classical multivariate Wick polynomials. In particular, we have that, in this case,
$$
  ev(W_T(a^{\otimes n}))=W_n(a)
$$
for a single element $a\in A$ \cite{EPTZ18}.
\end{remark}

Observe that by definition we have that ${\id}=(W_T\otimes\phi)\Delta^\shuffle$, so we get a tensor version of relation \eqref{eq:wickbasis}:
\begin{align*}
  a_1\dotsm a_n
  &= \sum_{S\subseteq[n]}W_T(a_S)\phi(a_{[n]\setminus S})\\
  &= \sum_{S\subseteq[n]}W_T(a_S)\sum_{\pi\in P([n]\setminus S)}\prod_{B\in\pi}c(a_B).
\end{align*}
Applying the evaluation map, the resulting relation is sometimes used as a recursive definition of the Wick polynomials \cite{HS17} in terms of moments or cumulants.

Since $\phi$ is not a character on $\overline{T}(A)$ (it is not multiplicative: $\phi(a_1a_2)=\phi(a_1 \cdot_{\!\scriptscriptstyle{A}} a_2)$ is different from the product $\phi(a_1)\phi(a_2)$ in general), it is not an element in the group of characters, i.e., the group-like elements in the completion of the dual graded Hopf algebra. Therefore, the $\exp/\log$ correspondence between tensor cumulants and moments cannot be analyzed from a Lie theoretic point of view. We refer the reader to \cite{EP17} for a discussion of the group and Lie algebra correspondence in the context of free probability.
The map $\phi$ has then a unique extension $\Phi\colon\overline T(T(A))\to \mathbb{K}$ as an algebra character. The unshuffle coproduct on the tensor algebra $\overline{T}(A)$ also admits a unique extension $\Delta^\shuffle\colon\overline T(T(A)) \to \overline T(T(A)) \otimes \overline T(T(A))$ as an algebra morphism:
$$
  \Delta^\shuffle(w_1|\dots |w_n)\coloneq\Delta^\shuffle(w_1)\cdot \dots \cdot\Delta^\shuffle(w_n),
$$
where the unit of $\overline{T}(A)$ is implicitly identified with the unit of $\overline{T}(T(A))$.

The following proposition and theorem are variants of the corresponding results in \cite{EPTZ18}, where they were obtained in the case where the algebra $A$ is commutative.

 \begin{proposition}\label{prop:antipode}
  The double tensor algebra $\overline T(T(A))$ with the coproduct $\Delta^\shuffle$ is a graded connected Hopf algebra, where $\deg(w_1\vert\dotsm\vert w_n)= \deg(w_1)+\dotsb+\deg(w_n)$. Its antipode $\mathcal S$ is the unique algebra anti-automorphism of $\overline T(T(A))$ such that
\[
  \mathcal S(a_1\cdots a_n)
    = \sum_{\pi\in P(n)}(-1)^{|\pi|}\sum_{\sigma\in\mathbb S_{|\pi|}}
    a_{B_{\sigma(1)}}\vert\dotsm\vert a_{B_{\sigma(|\pi|)}}
\]
for all $a_1,\dotsc,a_n\in A$.
\end{proposition}

As a consequence, we obtain using either Proposition \ref{prop:antipode} and lifting the computation of $W_T$ to $\overline T(T(A))$ or directly the definition of $W_T$:

\begin{theorem}
  The tensor Wick map admits the explicit expansion
\[
  W_T(a_1\dotsm a_n) = \sum_{S\subseteq[n]}a_S\sum_{\pi\in P([n]\setminus S)}(-1)^{|\pi|}|\pi|!\prod_{B\in\pi}\phi(a_B).
\]
\end{theorem}

Another point that is also addressed in \cite{EPTZ18} is the fact that Wick powers do not satisfy the usual rules of calculus: for example, since $\Wick{\!X\!}=X-\mathbb EX$ and $\Wick{\!X^2\!}=X^2-2X\mathbb EX+2(\mathbb EX)^2-\mathbb EX^2$ we see that
$\Wick{X}\!\!\cdot\Wick{X} \neq\Wick{X^2}\!\!$.
% Indeed, we have that
%\[
%  \Wick{X}\!\!\cdot\Wick{X}=X^2-2X\mathbb EX+(\mathbb EX)^2.
%\]
Nonetheless, using Hopf algebra techniques, the invertibility of the Wick map allowed us to define a modified product $\cdot_{\!\scriptscriptstyle{\varphi}}$ on polynomials such that
\[
  \Wick{X^n}\!\!\cdot_{\!\scriptscriptstyle{\varphi}}\Wick{X^m}=\Wick{X^{n+m}}\!\!,
\]
and a similar formula holds in the multivariate case. Since $W_T$ is a linear automorphism of
$\overline{T}(A)$,  these observations can be adapted to the tensor case as in
\Cref{sect:products}.
\bibliographystyle{arxiv}
\bibliography{../free_wick}
\end{document}